\documentclass[11pt,francais]{smfart}

\usepackage[T1]{fontenc}
\usepackage[english,francais]{babel}

\usepackage{amssymb,url,xspace,smfthm}

\newtheorem{e-proposition}[theorem]{Proposition}

\newtheorem{e-definition}[theorem]{Definition\rm}
\newtheorem{remark}{\it Remark\/}

\newtheorem{theoreme}{Th\'eor\`eme}[section]
\newtheorem{lemme}[theoreme]{Lemme}

\renewcommand{\theequation}{\arabic{equation}}
\newenvironment{preuve}[1][Preuve]{\noindent\textbf{#1.} }{\ \rule{0.5em}{0.5em}}
\begin{document}

\title{Sur une propri\'et\'e des polyn\^omes de stirling}
\date{22 F\'{e}vrier 2014}
\author{Farid BENCHERIF}
\address{Laboratoire LA3C, Facult\'e de Math\'emath\'ematiques, U.S.T.H.B., Bp 32 El Alia 16111 Bab Ezzouar Alger.}
\email{fbencherif@usthb.dz; fbencherif@yahoo.fr}
\author{Tarek GARICI}
\address{Laboratoire LA3C, Facult\'e de Math\'emath\'ematiques, U.S.T.H.B., Bp 32 El Alia 16111 Bab Ezzouar Alger.}
\email{tgarici@usthb.dz; tarekgarici@gmail.com}

\begin{abstract}
Dans cet article, nous r\'{e}pondons positivement \`{a} une question pos\'{e}%
e en 1960 par D.S. Mitrinovi\'{c} et R.S. Mitrinovi\'{c} \cite{mit6}
concernant les nombres de Stirling de premi\`{e}re esp\`{e}ce $s(n,k).$ Nous
prouvons que pour tout $k\geq 2$ il existe un entier $m_{k}$ et un polyn\^{o}%
me primitif $P_{k}\left( x\right) $ de $\mathbb{Z}\left[ x\right] $ tels que
pour tout $n\geq k$, $s\left(n,n-k\right)=\frac{1}{m_{k}}\binom{n}{k+1}\left(n(n-1)\right)^{\mathop{\rm mod}\nolimits(k,2)}P_{k}(n)$. De plus pour tout $k\geq1$, $P_{2k}(0)=P_{2k+1}(0)$.  
\end{abstract}

\begin{altabstract}
In this article, we give a positive answer to a question posed in 1960 by
D.S. Mitrinovi\'{c} and R.S. Mitrinovi\'{c} \cite{mit6} concerned the
Stirling numbers of the first kind $s(n,k).$ We prove that for all $k\geq 2$
there exist an integer $m_{k}$ and a primitive polynomial $P_{k}\left(
x\right) $ in $\mathbb{Z}[x]$ such that for all $n\geq k$, $s(n,n-k)=%
\frac{1}{m_{k}}\binom{n}{k+1}\left( n(n-1)\right) ^{\mathop{\rm mod}\nolimits%
(k,2)}P_{k}(n)$. Moreover for all $k\geq1$, $P_{2k}(0)=P_{2k+1}(0)$.
\end{altabstract}
\maketitle

\section{Introduction}

Le but de cet article est de r\'{e}pondre \`{a} une question pos\'{e}e par
D.S. Mitrinovi\'{c} et R.S. Mitrinovi\'{c} \cite{mit6} en 1960. Dans le Th%
\'{e}or\`{e}me qui suit,\ on a utilis\'{e} la notation $\left\lfloor
x\right\rfloor $ pour d\'{e}signer la partie enti\`{e}re d'un \ nombre r\'{e}%
el $x$. Un polyn\^{o}me non nul $a_{n}x^{n}+a_{n-1}x^{n-1}+\cdots
+a_{1}x+a_{0}$ de $%
\mathbb{Z}
\lbrack x]$ est dit primitif dans $\mathbb{Z}[x]$ si $\mathop{\rm pgcd}%
\nolimits(a_{k})_{0\leq k\leq n}=1$.

\begin{theoreme}
\label{theoreme}Soit $(m_{n})_{n\geq 0}$ la suite num\'{e}rique d\'{e}finie
par la relation 
\begin{equation}
m_{n}:=\frac{1}{(n+1)!}\prod\limits_{p\text{ premier et }p\text{ }\leqslant
n+1}p^{\left\lfloor \frac{n}{p-1}\right\rfloor +\left\lfloor \frac{n}{%
p\left( p-1\right) }\right\rfloor +\left\lfloor \frac{n}{p^{2}\left(
p-1\right) }\right\rfloor +\cdots }.  \label{th0}
\end{equation}%
Alors $(m_{n})_{n\geq 0}$ est une suite d'entiers naturels et, pour tout
entier $k\geqslant 1,$ on a 
\begin{equation}
s(n,n-2k)=\frac{1}{m_{2k}}\binom{n}{2k+1}P_{2k}(n)\text{ \ \ }(n\geqslant 2k)
\label{th1}
\end{equation}%
et 
\begin{equation}
s(n,n-2k-1)=\frac{1}{m_{2k+1}}\binom{n}{2k+2}n(n-1)P_{2k+1}(n)\text{ \ \ }%
(n\geqslant 2k+1),  \label{th2}
\end{equation}%
$P_{2k}(x)$ et\ $P_{2k+1}(x)$ \'{e}tant deux polyn\^{o}mes primitifs de $%
\mathbb{Z}
\lbrack x]$ v\'{e}rifiant la relation%
\begin{equation}
P_{2k}(0)=P_{2k+1}(0).  \label{th3}
\end{equation}
\end{theoreme}

La suite $\left( m_{n}\right) _{n\geq 0}=(1,1,4,2,48,16,576,144,...)$ d\'{e}%
finie en (\ref{th0}) est r\'{e}pertori\'{e}e sous la r\'{e}f\'{e}rence $%
A163176$ dans l'encyclop\'{e}die des suites d'entiers \cite{slo}. Pour $%
2\leqslant n\leqslant $ $9$, les expressions des polyn\^{o}mes $P_{n}(x)$
sont: 
\begin{eqnarray*}
P_{2}(x) &=&3x-1, \\
P_{3}(x) &=&-1, \\
P_{4}(x) &=&15x^{3}-30x^{2}+5x+2, \\
P_{5}(x) &=&-3x^{2}+7x+2, \\
P_{6}(x) &=&63x^{5}-315x^{4}+315x^{3}+91x^{2}-42x-16, \\
P_{7}(x) &=&-9x^{4}+54x^{3}-51x^{2}-58x-16, \\
P_{8}(x) &=&135x^{7}-1260x^{6}+3150x^{5}-840x^{4}-2345x^{3}-540x^{2} \\
&&+404x+144, \\
P_{9}(n) &=&-15x^{6}+165x^{5}-465x^{4}-17x^{3}+648x^{2}+548x+144.
\end{eqnarray*}%
D.S. Mitrinovi\'{c} et R.S. Mitrinovi\'{c} ont v\'{e}rifi\'{e} les relations
(\ref{th1}), (\ref{th2}) et (\ref{th3}) pour $k\in \left\{
1,2,3,4,5,6\right\} $. Ils ont alors propos\'{e} (\cite{mit6}, p. 4) le probl%
\`{e}me d'examiner si ces relations avaient lieu en g\'{e}n\'{e}ral pour
tout entier $k\geqslant 1$. Le Th\'{e}or\`{e}me r\'{e}pond positivement \`{a}
ce probl\`{e}me.

\section{D\'{e}monstration du Th\'{e}or\`{e}me \protect\ref{theoreme}}

La d\'{e}monstration du Th\'{e}or\`{e}me utilise trois lemmes et repose
essentiellement sur des propri\'{e}t\'{e}s des polyn\^{o}mes de N\"{o}rlund
et de la suite $\left( m_{n}\right) _{n\geq 0}$.

Les polyn\^{o}mes de N\"{o}rlund $B_{n}^{(x)}$ sont d\'{e}finis par (\cite%
{nor}, Chapitre 6)%
\begin{equation*}
\left( \frac{z}{e^{z}-1}\right) ^{x}=\sum_{n=0}^{\infty }B_{n}^{(x)}\frac{%
z^{n}}{n!}\text{.}
\end{equation*}%
$B_{n}^{(x)}$ est un polyn\^{o}me \`{a} coefficients rationnels de degr\'{e} 
$n$ divisible par $x$ pour $n\geqslant 1$. Les nombres de Bernoulli $B_{n}$
sont d\'{e}finis par $B_{n}=B_{n}^{(1)}$ \ \ $(n\geqslant 0).$On sait que%
\begin{equation}
B_{2n+1}=0\text{ \ \ \ }(n\geqslant 1).  \label{31}
\end{equation}%
\ 

\begin{lemme}
\label{lemm1}Pour $n\geqslant 1,$ on a%
\begin{equation}
\lbrack x]\left( B_{n}^{(x)}\right) =(-1)^{n-1}\frac{B_{n}}{n}  \label{rel19}
\end{equation}%
\begin{equation}
\lbrack x^{2}]\left( B_{2n+1}^{(x)}\right) =\frac{2n+1}{4n}B_{2n}
\label{rel22}
\end{equation}
\end{lemme}

\begin{preuve}
Soit $n\geqslant 1$. Dans (\cite{liu}, Th\'{e}or\`{e}mes $1$ et $2$), Liu et
Srivastava ont d\'{e}termin\'{e} explicitement\ les coefficients de $%
B_{n}^{(x)}$ en prouvant que le coefficient de $x^{k}$ dans $B_{n}^{(x)}$\
est donn\'{e} par%
\begin{equation}
\lbrack x^{k}]B_{n}^{(x)}=(-1)^{n-k}\frac{n!}{k!}\sum \frac{B_{\nu
_{1}}...B_{\nu _{k}}}{(\nu _{1}...\nu _{k})\nu _{1}!...\nu _{k}!}\text{ \ \ }%
(1\leqslant k\leqslant n),  \label{rel16}
\end{equation}%
la sommation ayant lieu sur les entiers $\nu _{1},...,\nu _{k}\geqslant 1$,
tels que $\nu _{1}+\cdots +\nu _{k}=n$. Soit $n\geqslant 1$. Pour $k=1$, (%
\ref{rel16}) permet d'obtenir ais\'{e}ment (\ref{rel19}). Pour $k=2$, (\ref%
{rel16}) permet d'\'{e}crire: 
\begin{eqnarray}
\lbrack x^{2}]\left( B_{2n+1}^{(x)}\right)  &=&-\frac{1}{2}\sum_{j=1}^{2n}%
\binom{2n+1}{j}\frac{B_{j}B_{2n+1-j}}{j(2n+1-j)}  \notag \\
&=&-\frac{1}{2}\binom{2n+1}{1}\frac{B_{1}B_{2n}}{2n}-\frac{1}{2}\binom{2n+1}{%
2n}\frac{B_{2n}B_{1}}{2n}  \notag \\
&&-\frac{1}{2}\sum_{j=2}^{2n-1}\binom{2n+1}{j}\frac{B_{j}B_{2n+1-j}}{%
j(2n+1-j)}  \notag \\
&=&\frac{2n+1}{4n}B_{2n}-\frac{1}{2}\sum_{j=2}^{2n-1}\binom{2n+1}{j}\frac{%
B_{j}B_{2n+1-j}}{j(2n+1-j)}.  \label{23}
\end{eqnarray}%
Ainsi la relation (\ref{rel22}) est bien v\'{e}rifi\'{e}e pour $n=1.$ Elle
l'est aussi pour $n\geqslant 2$ en remarquant que les termes figurant sous
le signe de sommation dans (\ref{23}) sont tous nuls car pour $n$ $\geqslant
2$ et $2\leqslant j\leqslant 2n-1,$ l'un au moins des deux nombres de
Bernoulli $B_{j}$ ou $B_{2n+1-j}$ est d'indice impair strictement plus grand
que $1$ et par suite $B_{j}B_{2n+1-j}=0$ d'apr\'{e}s (\ref{31}).
\end{preuve}

Le Lemme suivant est essentiel dans notre d\'{e}monstration.

\begin{lemme}
\label{lemm2}Pour $n\geqslant 2,$ on a 
\begin{equation}
\binom{x-1}{n}B_{n}^{(x)}=\frac{1}{m_{n}}\binom{x}{n+1}(x(x-1))^{\mathop{\rm
mod}\nolimits(n,2)}P_{n}(x)  \label{lem2}
\end{equation}%
o\`{u} $P_{n}(x)$ est polyn\^{o}me primitif de $%
\mathbb{Z}
\lbrack x]$.
\end{lemme}

\begin{preuve}
Soit $n\geqslant 1$. Pour tout nombre $p$ premier, d\'{e}signons par $%
r_{p}(n)$ l'exposant de la plus grande puissance de $p$ divisant $n!$.
Adelberg (\cite{ade1}, corollary 3) a montr\'{e} que si on pose: 
\begin{equation}
d_{n}=\frac{1}{n!}\prod\limits_{p\text{ }premier\text{ }et\text{ }p\leq
n+1}p^{r_{p}(n_{p})}\text{,}  \label{p1}
\end{equation}%
avec 
\begin{equation*}
n_{p}=p\left\lfloor \frac{n}{p-1}\right\rfloor ,
\end{equation*}%
alors 
\begin{equation*}
d_{n}B_{n}^{(x)}\text{est un polyn\^{o}me primitif de }%
\mathbb{Z}
\lbrack x].
\end{equation*}%
Par la formule de Legendre (\cite{ten}, p. 31), on a pour tout nombre
premier $p$ tel que $p\leqslant n+1$:%
\begin{equation}
r_{p}(n_{p})=\sum_{k\geq 0}\left\lfloor \frac{1}{p^{k}}\left\lfloor \frac{n}{%
p-1}\right\rfloor \right\rfloor =\sum_{k\geq 0}\left\lfloor \frac{n}{%
p^{k}(p-1)}\right\rfloor .  \label{p2}
\end{equation}%
De (\ref{p1}), (\ref{p2}) et (\ref{th0}), on d\'{e}duit que%
\begin{equation*}
d_{n}=(n+1)m_{n}.
\end{equation*}%
D'autre part, on sait que $B_{n}^{(x)}$ est divisible par $x$. Pour $n$
impair tel que $n\geqslant 3$, on a de plus, d'apr\`{e}s (\ref{rel19}) et (%
\ref{31}): 
\begin{equation*}
\lbrack x]\left( B_{n}^{(x)}\right) =(-1)^{n-1}\frac{B_{n}}{n}=0\text{ et }%
B_{n}^{(1)}=B_{n}=0.
\end{equation*}%
Il en r\'{e}sulte que dans $%
\mathbb{Z}
\lbrack x]$, le polyn\^{o}me primitif $(n+1)m_{n}B_{n}^{(x)}$est divisible
par le polyn\^{o}me primitif $x\left( x(x-1)\right) ^{\mathop{\rm mod}%
\nolimits(n,2)}$ pour $n\geqslant 2.$ Le quotient $P_{n}(x)$ de ces deux
polyn\^{o}mes est aussi un \ polyn\^{o}me primitif de $%
\mathbb{Z}
\lbrack x]$ et on a donc: 
\begin{equation}
(n+1)m_{n}B_{n}^{(x)}=x\left( x(x-1)\right) ^{\mathop{\rm mod}\nolimits%
(n,2)}P_{n}(x)\text{ \ \ }(n\geqslant 2).\text{ }  \label{r1}
\end{equation}%
En multipliant les deux membres de (\ref{r1}) par $\frac{1}{(n+1)m_{n}}%
\binom{x-1}{n},$ on obtient (\ref{lem2}).
\end{preuve}

Le Lemme suivant \'{e}nonce des propri\'{e}tes de la suite num\'{e}rique $%
\left( m_{n}\right) _{n\geq 0}$ \ d\'{e}finie en (\ref{th0}).

\begin{lemme}
\label{lemm3}Pour tout entier $n\geq 0$

\begin{enumerate}
\item $m_{n}$ est un entier

\item $m_{2n}=(n+1)m_{2n+1}$
\end{enumerate}
\end{lemme}

\begin{preuve}
Soit $n\geq 0$ un entier.

\begin{enumerate}
\item Pour tout nombre premier $p\leqslant n+1$ et pour tout entier $%
k\geqslant 0$, on a%
\begin{equation*}
\frac{n}{p^{k}(p-1)}-\frac{n+1}{p^{k+1}}=\frac{n+1-p}{p^{k+1}(p-1)}\geqslant
0
\end{equation*}%
et par cons\'{e}quent%
\begin{equation*}
\left\lfloor \frac{n}{p^{k}(p-1)}\right\rfloor -\left\lfloor \frac{n+1}{%
p^{k+1}}\right\rfloor \geqslant 0.
\end{equation*}%
Il en r\'{e}sulte que $m_{n}$ est un entier. En effet, par la formule de
Legendre, on a%
\begin{equation*}
v_{p}\left( m_{n}\right) =\sum_{k\geq 0}\left( \left\lfloor \frac{n}{%
p^{k}(p-1)}\right\rfloor -\left\lfloor \frac{n+1}{p^{k+1}}\right\rfloor
\right) \geqslant 0.
\end{equation*}

\item Soit un nombre premier $p\leqslant n+1$. On montre ais\'{e}ment que
pour tous entiers naturels non nuls $x$ et $y:$%
\begin{equation*}
\left\lfloor \frac{x+1}{y}\right\rfloor -\left\lfloor \frac{x}{y}%
\right\rfloor =\left\{ 
\begin{array}{l}
1\text{ \ si }y\text{ divise }x+1, \\ 
0\text{ \ sinon.}%
\end{array}%
\right.
\end{equation*}%
Il en r\'{e}sulte que si $p=2,$ on a 
\begin{eqnarray*}
v_{2}\left( \frac{(n+1)m_{2n+1}}{m_{2n}}\right) &=&\sum_{k\geq
1}\left\lfloor \frac{2n+1}{2^{k}}\right\rfloor -\left\lfloor \frac{2n}{2^{k}}%
\right\rfloor \\
&=&0
\end{eqnarray*}%
et si $p\geqslant 3,$ on a aussi 
\begin{eqnarray*}
v_{p}\left( \frac{(n+1)m_{2n+1}}{m_{2n}}\right) &=&\sum_{k\geq
0}\left\lfloor \frac{2n+1}{p^{k}(p-1)}\right\rfloor -\left\lfloor \frac{2n}{%
p^{k}(p-1)}\right\rfloor \\
&=&0,
\end{eqnarray*}%
car $p^{k}(p-1)$ est alors un entier pair et il ne peut donc pas diviser $%
2n+1$. Par suite on a pour tout nombre premier $p,$ $v_{p}\left( \frac{%
(n+1)m_{2n+1}}{m_{2n}}\right) =0$, ce qui \'{e}quivaut \`{a} affirmer que%
\begin{equation*}
\frac{(n+1)m_{2n+1}}{m_{2n}}=1.
\end{equation*}
\end{enumerate}
\end{preuve}

Nous pouvons maintenant prouver le Th\'{e}or\`{e}me. Soit $k$ un entier sp%
\`{e}rieure \`{a} $1$. Il est bien connu que l'on a (\cite{char}, p. 329): 
\begin{equation}
s(n,n-j)=\binom{n-1}{j}B_{j}^{(n)},\text{ pour }n\geq j\geq 0.  \label{150}
\end{equation}%
Avec le Lemme \ref{lemm2}, (\ref{150}) s'\'{e}crit 
\begin{equation}
s(n,n-j)=\frac{1}{m_{j}}\binom{n}{j+1}(n(n-1))^{\mathop{\rm mod}\nolimits%
(j,2)}P_{j}(n)\text{ \ \ }(n\geq j\geq 2),  \label{151}
\end{equation}%
$P_{j}(x)$ \'{e}tant un polyn\^{o}me primitif de $%
\mathbb{Z}
\lbrack x].$ Pour $j=2k$ (resp $j=2k+1)$, la relation (\ref{151}) se traduit
par (\ref{th1}), (resp (\ref{th2})).

De plus, en choisissant $n=2k$ puis $n=2k+1$ dans la relation (\ref{r1}), on
obtient%
\begin{equation*}
\left\{ 
\begin{array}{l}
(2k+1)m_{2k}B_{2k}^{(x)}=xP_{2k}(x), \\ 
(2k+2)m_{2k+1}B_{2k+1}^{(x)}=x^{3}P_{2k+1}(x)-x^{2}P_{2k+1}(x)\text{.}%
\end{array}%
\right.
\end{equation*}%
On en d\'{e}duit que%
\begin{equation*}
\left\{ 
\begin{array}{l}
P_{2k}(0)=[x]((2k+1)m_{2k}B_{2k}^{(x)}), \\ 
P_{2k+1}(0)=[x^{2}](-(2k+2)m_{2k+1}B_{2k+1}^{(x)})\text{.}%
\end{array}%
\right.
\end{equation*}%
A la lumi\`{e}re du Lemme \ref{lemm1}, ces deux derni\`{e}res relations
deviennent%
\begin{equation}
\left\{ 
\begin{array}{l}
P_{2k}(0)=-m_{2k}(2k+1)\frac{B_{2k}}{2k}, \\ 
P_{2k+1}(0)=-(k+1)m_{2k+1}(2k+1)\frac{B_{2k}}{2k}.%
\end{array}%
\right.  \label{200}
\end{equation}%
La relation 2 du Lemme \ref{lemm3} permet alors de d\'{e}duire de (\ref{200}%
) que \newline
$P_{2k+1}(0)=P_{2k}(0)$, ce qui \'{e}tablit (\ref{th3}). La d\'{e}%
monstration du Th\'{e}or\`{e}me est compl\`{e}te.

\begin{remark}
Comme on sait que le degr\'{e} de $B_{n}^{(x)}$est \'{e}gal \`{a} $n$ pour
tout $n\geqslant 0$, on en d\'{e}duit \`{a} l'aide de la relation (\ref{r1})
que pour $k\geqslant 1$, on a:%
\begin{equation*}
\deg (P_{2k})=2k-1\text{ et }\deg (P_{2k+1})=2k-2.
\end{equation*}
\end{remark}

\end{document}